\documentclass[twoside]{amsart}
\usepackage{ulem} 
\usepackage{latexsym}
\usepackage{amssymb,amsmath,amsopn}
\usepackage[dvips]{graphicx}   
\usepackage{color,epsfig}      
\usepackage{enumitem}
\usepackage{caption}
\usepackage{fourier}
\usepackage[cp1251]{inputenc}
               {\begin{list}{}{\leftmargin#1\rightmargin#2}\item{}}%
               {\end{list}}

\newtheorem{thm}{Theorem}
\newtheorem{lem}{Lemma}
\newtheorem{prop}{Proposition}
\theoremstyle{definition}
\newtheorem{defn}{Definition}
\newtheorem*{rem}{Remark}

\newtheorem{tcor}{Corollary}[thm]

\renewcommand{\P}{\mathcal{P}}
\newcommand{\Q}{\mathcal{Q}}

\newcommand{\Tet}{\mathcal T}
\newcommand{\Pol}{\mathcal P}
\newcommand{\scrA}{\mathcal{A}}
\newcommand{\scrB}{\mathcal{B}}
\newcommand{\scrC}{\mathcal{C}}
\newcommand{\scrD}{\mathcal{D}}

\definecolor{green}{rgb}{0.0, 0.75, 0.25} 

\def\eea{\end{eqnarray}}
\renewcommand{\emph}{\textit}

\DeclareMathOperator{\aff}{aff}
\DeclareMathOperator{\relint}{relint}

\DeclareMathOperator{\inter}{int}
\DeclareMathOperator{\bd}{bd}

\title[Equilibria of tetrahedra]{On equilibria of tetrahedra}
\author[G. Alm\'adi, R. J. MacG. Dawson, G. Domokos and K. Reg\H{o}s ] {Gerg\H o Alm\'adi, Robert J. MacG. Dawson, G\'abor Domokos and Krisztina Reg\H{o}s}
\address{Gerg\H o Alm\'adi, ELKH-BME Morphodynamics Research Group, Budapest University of Technology and Economics,
M\H uegyetem rakpart 1-3., Budapest, Hungary, 1111}
\email{gergo.almadi14@gmail.com}
\address{Robert J. MacG. Dawson, Dept. Mathematics and Computer Science, Saint Mary's University, Halifax, Nova Scotia B3H 3C3, Canada}
\email{rdawson@cs.stmarys.ca}
\address{G\'abor Domokos, Dept. of Morphology and Geometric Modeling and ELKH-BME Morphodynamics Research Group, Budapest University of Technology and Economics,
M\H uegyetem rakpart 1-3., Budapest, Hungary, 1111}
\email{domokos@iit.bme.hu}
\address{Krisztina Reg\H os, Dept. of Morphology and Geometric Modeling and  ELKH-BME Morphodynamics Research Group,
M\H uegyetem rakpart 1-3., Budapest, Hungary, 1111}
\email{regos.kriszti@gmail.com}

\thanks{GD and GA: Support of the NKFIH Hungarian Research Fund grant 134199 and of grant BME FIKP-V\'IZ by EMMI is kindly acknowledged. KR:  This research has been supported by the
program UNKP-22-3 by ITM and NKFIH. The gift representing the Albrecht Science Fellowship
is gratefully appreciated.}
\subjclass[2010]{52B10, 77C20, 52A38}

\keywords{polyhedron, static equilibrium, monostatic polyhedron}

\begin{document}
\maketitle

\begin{abstract}
The monostatic property of polyhedra (i.e. the property of having just one stable or unstable static equilibrium point) has been in a focus of research ever since Conway and Guy \cite{Conway} published the proof of the existence of the first such object. In the same article they also proved  that a homogeneous tetrahedron has at least two stable equilibrium points.  By using polar duality,  the same idea has been used \cite{balancing} to prove that a homogeneous tetrahedron has at least two unstable equilibria. Conway \cite{Dawson} also claimed that among inhomogeneous tetrahedra one can find monostable ones.  Here we not only give a formal proof of this statement and show that monostatic tetrahedra have exactly 4 equilibria, but also demonstrate a startling new aspect of this problem: being monostatic implies certain \emph{visible} features of the shape and vice versa.  Our results also imply that mono-monostatic tetrahedra (having just one stable and just one unstable equilibrium point) do not exist.  In contrast, we show that for any other legal number of faces, edges, and vertices there is a mono-monostatic polyhedron with that face vector.
\end{abstract}


\section{Introduction}

\subsection{History and background}
 The curious mechanical properties of tetrahedral solids were first studied in 1967 when Heppes constructed a homogeneous tetrahedron \cite{Heppes1967} that could rest stably on only two of its faces. Conway showed \cite{Dawson} that for a homogeneous tetrahedron this number cannot be improved: i.e., it cannot be \emph{monostable}.
 
 A three-dimensional, weighted convex polyhedron $\P$ with center of mass $O$, supported by a fixed horizontal plane, has three sorts of equilibria: stable (on a face), unstable (on a vertex), and saddle (on an edge). These correspond to local minima, maxima, and saddle points of the radial height function  $r_{\P,O}:S^2\rightarrow \mathbb{R}^+$ describing the boundary of $\P$ as a distance measured from $O$. The global study of such equilibrium points, today associated with Morse theory, goes back (on smooth surfaces interpreted in a Cartesian coordinate system) to Cayley \cite{Cayley} in 1859. Maxwell \cite{Maxwell} noted a few years later that Euler's formula $f-e+v=2$  described the relationship between maxima, saddles, and minima;  applied to the radial height function $r_{\P,O}$, this links the numbers of those faces, edges, and vertices of the convex polyhedron ${\P}$ which carry equilibria.  
 
 Conway claimed that monostability is possible for weighted tetrahedra, and asked whether it was possible for homogeneous simplices in higher dimensions. This was answered in a series of papers  \cite{Dawson, Dawson2, Dawson3} by Dawson and Finbow, who also showed \cite{Loaded} that even some regular polytopes, appropriately weighted, can be monostable. 

Here we focus on the weighted case. We will exhibit conditions on a (non-regular) tetrahedron $\Tet$ equivalent to the existence of a weighting $(\Tet,O)$ making it monostable; and we will offer necessary and sufficient conditions for a weighted tetrahedron to illustrate Conway's claim. (While this is not hard to construct, there does not appear to be a published example.) 

Although never stated explicitly, one underlying question about monostability has always been whether it could be a \emph{visible} property of any given shape. In general,  this  does not seem to be the case: the extreme (0.1\%) shape tolerance of the G\"omb\"c shape signal that it would be hard to select a mono-monostatic shape based alone on visual inspection. Indeed, in case of smooth shapes it was shown \cite{VarkonyiDomokos} that mono-monostatic ones have minimal flatness and thinness, thus they appear in the vicinity of the sphere. Here we show a  startling feature of tetrahedra, where both monostability and mono-instability are reflected in the shape in a unique manner. Not only are these extrinsic features beautiful, they also appear to be at the heart of the phenomenon: building on these shape characterictics
we will show that any monostable tetrahedron is \emph{bi-unstable} - that is, it has equilibrium on exactly two vertices. Neither G\"{o}mb\"{o}c-like mono-monostatic weighted tetrahedra, nor monostable tetrahedra with three or four unstable equilibria, exist.  However, we will show that if a tetrahedron has a weighting making it monostable on one face, then it must have other weightings making it monostable on each of its faces. We will also show that while, in general, the physics of tipping bodies in three or more dimensions is highly complicated, the tipping of a monostable tetrahedron can be completely described knowing its shape and center of mass.

The \emph{polar dual} of a polyhedron $\P$ is the set $\P^* = \{x:x\cdot p \leq 1 \mbox{ for all } p\in \P\}$. If $\P$ is convex,  $\P^{**}=\P$. The polar dual of a tetrahedron is also a tetrahedron; and there is a natural pairing between the vertices of one and the faces of the other. The following result was proved in \cite{balancing}.

\begin{prop}\label{PolarDualProp}
Let $(\P,O)$ be a weighted convex polyhedron, with $O$ at the origin. Then $\P$ has an equilibrium on a face if and only if $\P^*$ has an equilibrium on the corresponding vertex.
\end{prop}

Thus, any mono-unstable weighted tetrahedron is bistable. We will exhibit explicit geometric conditions for a tetrahedron to have such a weighting; and we can see that (in contrast with the monostable case) if $(\Tet,O)$ is mono-unstable on a vertex $A$, $\Tet$ cannot be weighted to be mono-unstable on any other vertex.  Finally, we will show that the face vector $(f,e,v) = (4,6,4)$ that characterizes tetrahedra is unique among those of polyhedra, in that any other legal face vector does have a representative polyhedron that may be weighted to make it mono-monostatic.

\subsection{Definitions of equilibrium}

\begin{defn}\label{def1}\cite{balancing}, \cite{VarkonyiDomokos}
Let $\Pol$ be a convex polyhedron, and let $\inter \Pol$ and $\bd \Pol$ denote its interior and boundary, respectively. We select a point $O \in \inter \Pol$, which we shall think of as the center of mass. (We are not assuming $\Pol$ to have uniform density, so this implies no restriction on $O$ other than its being an interior point.) 

We say that $(\Pol,O)$ is \emph{in equilibrium on} a face, edge, or vertex $A$ if there exists $Q \in \relint A$ such that the plane perpendicular to $[O,Q]$ at $Q$ supports $\Pol$.  (Recall that ``relative interior'' is defined in such a way that a singleton's relative interior is itself, though its \emph{interior} is empty: thus $\Pol$ may be in equilibrium on a vertex.)  We call the equilibrium \emph{stable} if $Q$ is on the relative interior of a face, \emph{unstable} if $Q$ is a vertex, and \emph{hyperbolic} (saddle) otherwise and we denote their numbers by $S,U,H$, respectively. 
\end{defn}

As noted above, Maxwell showed that $S - H + U = 2$. These equilibria correspond intuitively to positions in which a physical model of $(\Pol,O)$ balances on a horizontal surface.  They also  correspond to ``pits'', ``peaks'', and ``passes'' in the radial function of $\Pol$ with respect to $O$.

\begin{defn}\label{def2}\cite{balancing}
We call a $\Pol$ convex polyhedron \emph{monostable} if it has a unique stable equilibrium (there is exactly one face upon which it will rest) and \emph{mono-unstable} if it has a unique unstable equilibrium (there is exactly one vertex upon which it can balance precariously.) 
\end{defn}

\section{Results on monostability}

Henceforth we assume that $\Pol$ is a tetrahedron $\Tet = \boxtimes ABCD$ with face $\scrA$ opposite vertex $A$ (etc.) We say that a tetrahedron has an \emph{obtuse path} $A-B-C-D$ if it has three edges, WLOG $\overline{AB}$,$\overline{BC}$, and  $\overline{CD}$, with obtuse dihedral angles and no common vertex. Such tetrahedra exist, for instance tetrahedron $\Tet _0$ with vertices  $(A,B,C,D) = ((0,0,0),(0,0,100000),(153600,44400,0),(112200,7800,6400)$. (We remark
that $\Tet _0$ is monostable on face $\scrD$ with the center of mass located at
$O = (104200,4300,100).$)
We note that three obtuse edges can never surround a face, and that no tetrahedron can have four obtuse dihedrals.

\begin{thm}\label{MonoThm}
Let $\Tet$ be a tetrahedron; then the following are equivalent:
     \begin{enumerate}
         \item $\Tet$ has an obtuse path;
         \item there exists $O$ such that $(\Tet,O)$ is monostable;
         \item for every face $F$, there exists $O_F$ such that $(\Tet,O_F)$ is monostable on $F$.
     \end{enumerate}
\end{thm}
\begin{proof}
(3)$\Rightarrow$(2) trivially. 

(2) $\Rightarrow$ (1): A tetrahedron has the property (unique among polyhedra) that we can walk from vertex to vertex along a path of edges if and only if we can skip from face to face across the same edges. There must be enough obtuse dihedrals to let the tetrahedron roll from any face to the resting face; and no face can have three obtuse dihedrals. The obtuse dihedrals thus form a path of length 3.

Finally, (1)$\Rightarrow$(3). If our obtuse path is $A-B-C-D$, the obtuse edges connect the faces in the order $\scrC-\scrD-\scrA-\scrB$ (Figure \ref{Monostable}). It suffices to show that for appropriate $O_{\scrA}$ the pair $(\Tet,O_{\scrA})$ is monostable on $\scrA$, and similarly for $\scrB$.  Construct the plane perpendicular to $\scrC = \triangle ABD$ along the edge $\overline{AB}$ shared with $\scrD$. This cuts $\overline{CD}$ at a point $E$ (Figure \ref{Monostable}a). The tetrahedron $\boxtimes ABCE$  has obtuse edges $\overline{BC}$ and $\overline{EC}$. If $O \in \inter\boxtimes ABCE$, then $(\Tet,O)$ has no stable equilibrium on $\scrC$. Next, construct the plane perpendicular to $\scrD = \triangle ABC$ along $\overline{BC}$. This cuts $\overline{AE}$ at $F$, and $\overline{CE}$ is an obtuse edge of the tetrahedron $\boxtimes BCEF$ (Figure \ref{Monostable}b).  If $O \in \inter\boxtimes BCEF$, then $(\Tet,O)$ has no stable equilibrium on $\scrC$ or $\scrD$.

\begin{figure}[h!]
\centering
    \includegraphics[width=0.95\textwidth]{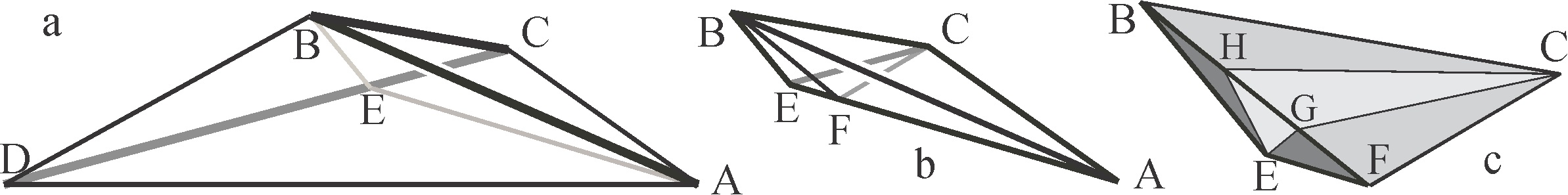}
\caption{Loading regions for a tetrahedron with an obtuse path}
\label{Monostable}
\end{figure}

If we repeat this with a plane perpendicular to $\scrA$ along $\overline{CE}$, it meets $\overline{BF}$ at $G$; and if $O_\scrB \in \inter\boxtimes CEFG$, then $(\Tet,O_\scrB)$ has no stable equilibrium on $\scrA$, $\scrC$, or $\scrD$, so is monostable on $\scrB$ (Figure \ref{Monostable}c).  Similarly, the plane through the same edge but perpendicular to $\scrB$ meets $\overline{BF}$ at $H$; and if $O_\scrA \in \inter\boxtimes BCEH$, then $(\Tet,O_\scrA)$ has no stable equilibrium on $\scrB$, $\scrC$, or $\scrD$ and is monostable on $\scrA$.
\end{proof}
\begin{rem} We've shown that, suitably weighted, a tetrahedron with an obtuse path is stable only on one face. We haven't shown how it gets there. (As the bartender says at closing time, "You don't have to go home, ladies and gentlemen, but you can't stay here.") In fact, without some way to dissipate energy, the tetrahedron will never settle onto any face, but will bounce forever. \footnote{Fans of opera, or at least of operatic trivia, will recall the story of Eva Turner, in the role of Tosca, throwing herself from the battlements onto an over-resilient trampoline placed there by the stage hands, and making several unplanned curtain calls!} 

When tipping over an edge, a polyhedron has only one degree of freedom, and how it dissipates its energy does not affect where it ends up, so long as it does so effectively.  However, in some cases, the body may tip not onto an edge (from where it must continue to the next face) but onto a vertex (figure \ref{FigTip} a). Should this happen, the body temporarily has not two but three degrees of freedom: the center of mass $O$ moves on a sphere about $A$, and the body can also rotate about the axis $OA$ (figure \ref{FigTip} b). We thus need to take into account torque, moment of inertia, and the position of $O$ relative to the  edge upon which the tetrahedron lands: the problem would seem intractable.

\begin{figure}[h!]
\centering
    \includegraphics[width=0.65\textwidth]{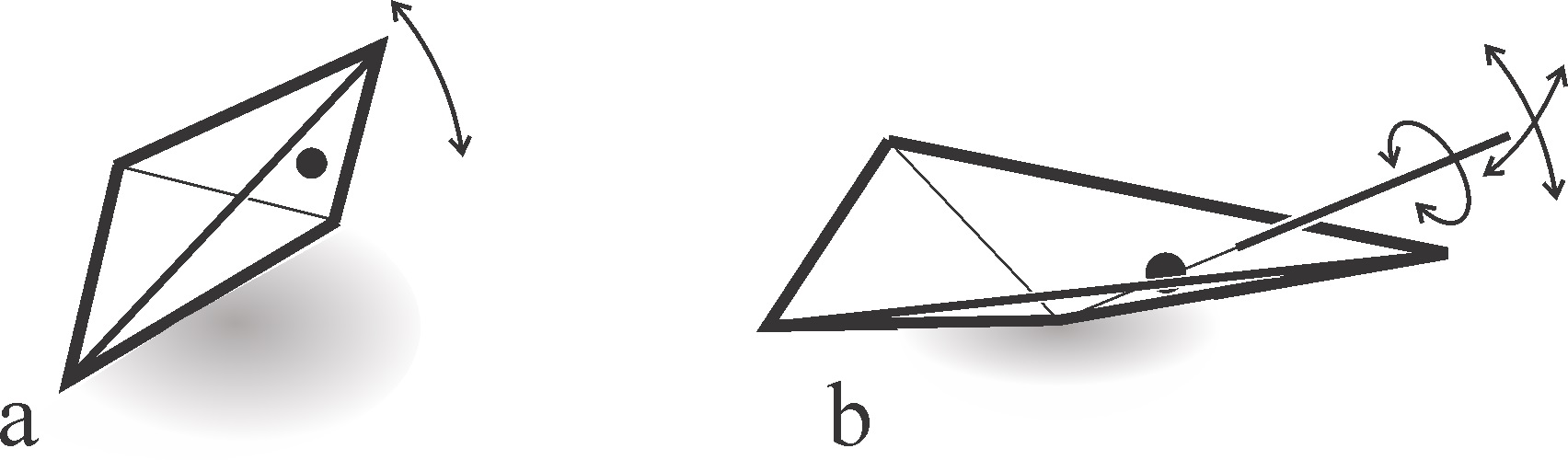}
\caption{A tetrahedron that rolls without slipping can have 1 or 3 degrees of freedom.}
\label{FigTip}
\end{figure}

Fortunately, these difficulties never arise if $(\Tet,O)$ is a monostably-weighted tetrahedron! In this case, as shown above, if a face has two obtuse dihedrals, the center of mass is positioned so that $\Tet$ will tip \emph{onto} that face across at least one of those edges. Each non-equilibrium face thus has a unique exit; and providing that we assume landings to be inelastic, the exact path to stable equilibrium may be computed knowing only $(\Tet,O)$ and the starting face.

\end{rem}

\section{Some spherical geometry}

\begin{figure}[h!]
\centering
    \includegraphics[width=0.4\textwidth]{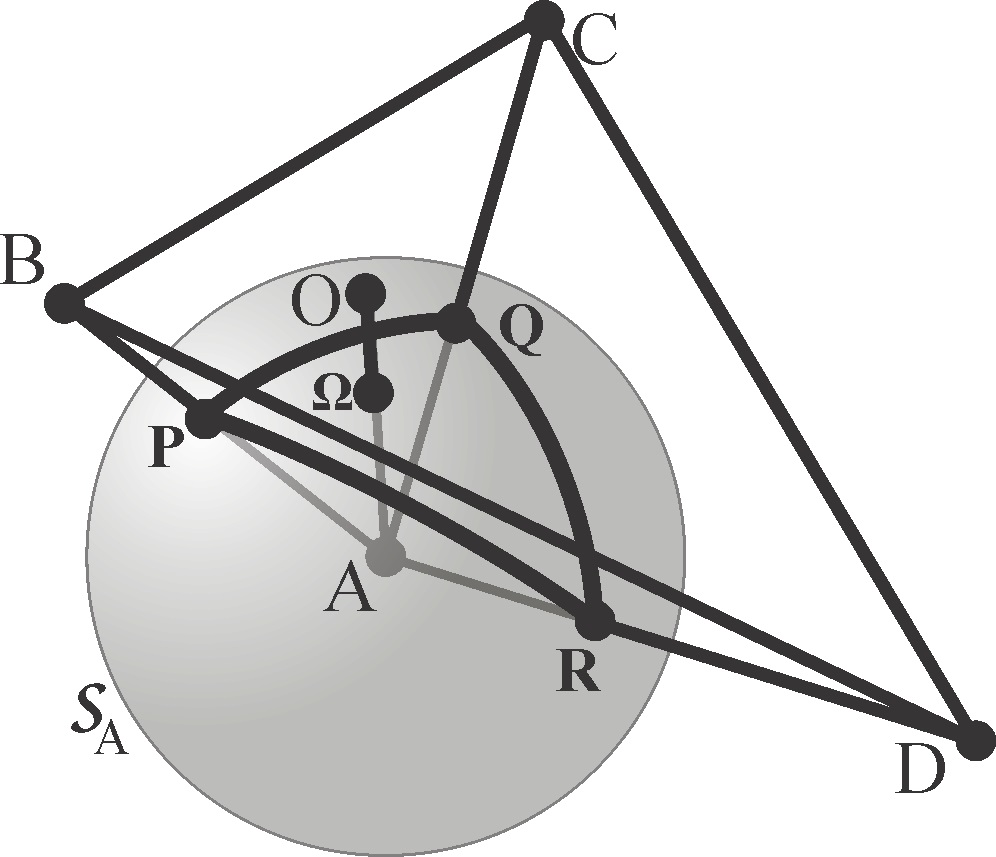}
\caption{The geometry of a vertex is the geometry of a sphere.}
\label{SphericalVertex}
\end{figure}

If we consider the intersection of $\Tet$ with a small sphere $\mathcal{S}_A$ centered at some vertex $A$, we see that the geometry of polyhedral vertices is just that of the sphere! For a weighted tetrahedron $(\Tet,O)$, let $P,Q,R,\Omega$ be the respective intersections of $\mathcal{S}_A$ with $\overline{AB},\overline{AC},\overline{AD}$, and $\overline{AO}$ (see Figure \ref{SphericalVertex}.) Then (for instance) the face angle $\angle BAC$ corresponds to the arc $\overline{PQ}$ on $\mathcal{S}_A$, while the dihedral angle between  faces $\triangle ABC$ and $\triangle ACD$ corresponds to the spherical angle $\angle PQR$. In each case the radian measures are equal. We call a spherical segment \emph{short} if its measure is less than $\pi/2$, otherwise \emph{long}; angles, as usual, are ``acute'' or ``obtuse.''

The following result characterizes unstable equilibria.

\begin{lem}\label{UnstableCond}
For any vertex $A$ of a weighted tetrahedron $(\Tet,O)$
\begin{enumerate}
    \item $(\Tet,O)$ has (unstable) equilibrium on  $A$  if and only if the angles $\angle BAO$,  $\angle CAO$, and  $\angle DAO$ are all acute. Equivalently, all of the spherical arcs $\overline{P\Omega}$, $\overline{Q\Omega}$, or $\overline{R\Omega}$ are short.
    \item  $(\Tet,O)$ has an equilibrium on $A$ for every $O\in\inter\Tet$ if none of the face angles $\angle BAC$,  $\angle CAD$,  $\angle DAB$ are obtuse. Equivalently, none of the arcs $\overline{PQ}$, $\overline{QR}$, or $\overline{RP}$ is long.  
\end{enumerate}
\end{lem}
\begin{proof}  Let $\Pi$ be the plane normal to $\overline{OA}$ at $A$: then $\Tet$ has an equilibrium on $A$ if and only if $B$,$C$, $D$ all lie on the same side of $\Pi$ as $O$. This is true for every $O$ \emph{interior to} $\Tet$  if and only if $\angle BAC$,  $\angle CAD$,  $\angle DAB$ are all acute or right.
\end{proof}

We can also characterize stable equilibria in this way, though we have local configurations at three vertices to consider. The following result is obvious:

\begin{lem}
\begin{enumerate}
    \item $(\Tet,O)$ has stable equilibrium on $\triangle ABC$ if and only if  the dihedral angles between (on the one hand) $\triangle ABC$ and (on the other hand) $\triangle ABO$, $\triangle AOC$, and $\triangle OBC$ are all acute.
    \item  $(\Tet,O)$ has stable equilibrium on $\triangle ABC$ for every $O\in\inter\Tet$ if and only if none of the dihedral angles between (on the one hand) $\triangle ABC$ and (on the other hand) $\triangle ABD$, $\triangle ADC$, or $\triangle DBC$ are obtuse.
    \item If (for instance) the dihedral angle between $\triangle ABC$ and  $\triangle ABO$ is obtuse, then (on $\mathcal{S}_A$) the angle $\angle QP\Omega$ is obtuse, as is the corresponding angle on the sphere $\mathcal{S}_B$.
\end{enumerate}
\end{lem}

Now that we've seen the significance of spherical geometry to this problem, let's establish a few facts from the folklore.

\begin{lem}\label{SmallProp} 
\begin{enumerate}
    \item A spherical triangle with only acute angles has only short edges; 
    \item A spherical triangle with exactly one acute angle has exactly one short edge, which is opposite the acute angle;
    \item  A spherical triangle with three long edges has three obtuse angles;
    \item A spherical triangle with exactly one long edge has exactly one obtuse angle, opposite the long edge.
    \item A spherical triangle with only short edges has at most one obtuse angle.
\end{enumerate}
 \end{lem}
\begin{proof}
\begin{enumerate}
\item If $\triangle ABC$ has only acute angles, $\cos A$, $\cos B$, and $\cos C$ are all positive. Then if (for instance) edge $BC$ has radian length $a$, one of the spherical cosine laws gives us that 
$$ \cos a = \frac{\cos A  + \cos B \cos C}{\sin B \sin C} > 0,$$
whence $a < \pi/2$. The proofs for $b$ and $c$ are similar.
\item If $A$ is acute, let $A'$ be the antipodal point: the colunar triangle $\triangle A'BC$ satisfies the conditions of (1).
\item If $\triangle ABC$ has only long edges, $\cos a$, $\cos b$, and $\cos c$ are all negative; the other spherical cosine law gives
$$ \cos A = \frac{\cos a - \cos b \cos c}{\sin b \sin c} < 0$$
and $A$ is obtuse: the proofs for $B$ and $C$ are similar.
\item again follows from (3) by consideration of the colunar triangle.
\item In this case, $\cos a$, $\cos b$, and $\cos c$ are all positive. If $A$ is obtuse, $\cos a < \cos b \cos c$, whence $a$, opposite $A$, must be the strictly longest edge.
\end{enumerate}
\end{proof}

We can, however, construct spherical triangles with exactly one obtuse angle and zero, or two, long edges. We can also construct a spherical triangle with three obtuse angles and only two long edges.  These results are summarized in Table \ref{AngleTable}. 
\begin{table}
\begin{center}
\begin{tabular}{|l|c|c|c|c|c|}
\hline
&\multicolumn{5}{|c|}{Obtuse angles}\\
&\multicolumn{5}{|c|}{(dihedrals)}\\
\hline
Long&&0&1&2&3 \\
\cline{2-6}
edges&0&$\sqrt{}$&$\sqrt{}$&x&x\\
\cline{2-6}
(obtuse&1&x&$\sqrt{}$&x&x\\
\cline{2-6}
face&2&x&$\sqrt{}$&$\sqrt{}$&$\sqrt{}$\\
\cline{2-6}
angles)&3&x&x&x&$\sqrt{}$\\
\hline
\end{tabular}\caption{Possible combinations of long/obtuse elements in spherical triangles and polyhedral vertices}\label{AngleTable}
\end{center}\end{table}

\section{Results on instability}

We begin by ruling out the possibility of a ``weighted tetrahedral G\"{o}mb\"{o}c'', and in fact prove more.
\begin{thm}\label{NoMonoMono}
No tetrahedron has both a monostable weighting and a mono-unstable weighting, even with different centers of mass.
\end{thm}

\begin{proof}
    As observed above, every monostable weighted tetrahedron $(\Tet,O)$ has two vertices $B,C$ that each have two obtuse dihedrals.  Table \ref{AngleTable} shows that each of these two must have two obtuse face angles; but a tetrahedron cannot have more than four obtuse face angles in total, so the other two vertices $A,D$ have only acute face angles. By Lemma \ref{UnstableCond}, $(\Tet,O')$ has equilibria on those vertices for any $O' \in \inter\Tet$.   
\end{proof}

\begin{figure}[h!]
\centering
    \includegraphics[width=0.95\textwidth]{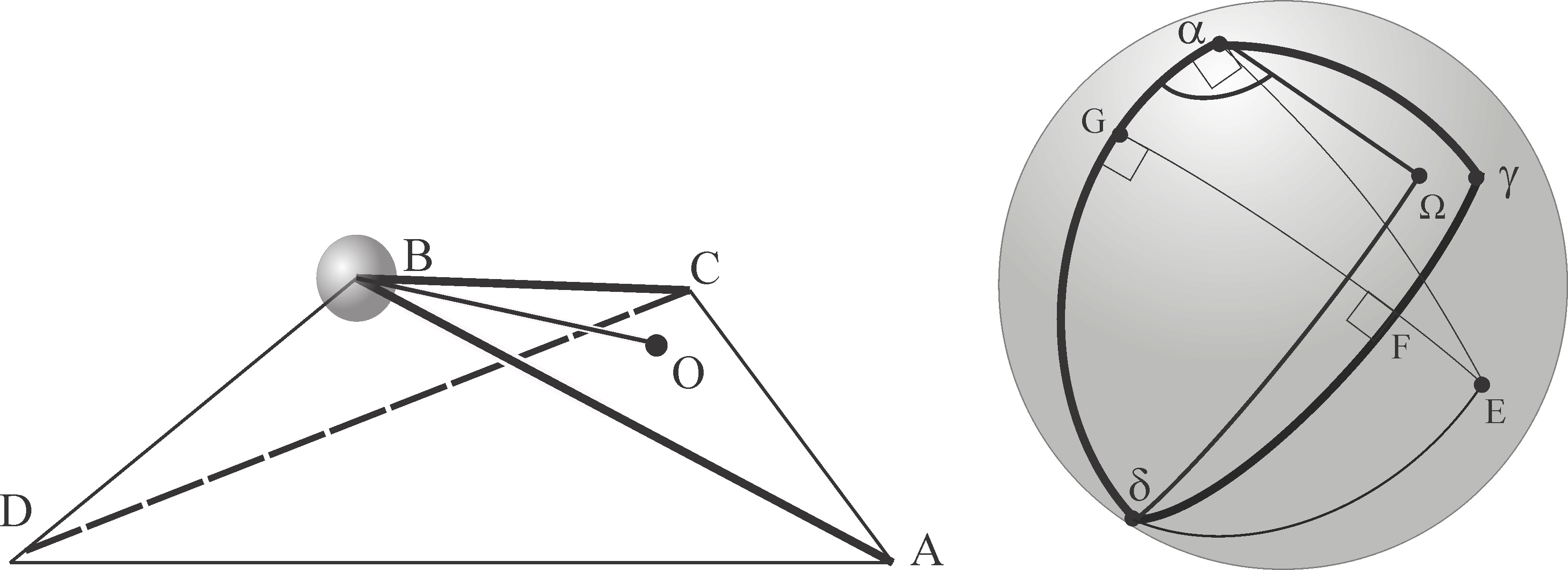}
\caption{The configuration at a vertex with two obtuse dihedral angles.}
\label{MonoThmFig}
\end{figure}

For a specific weighting, we can say more: 
\begin{thm}\label{MonoBi}
  If a weighted tetrahedron $(\Tet,O)$ is monostable, it has unstable equilibrium on exactly two vertices.
\end{thm} 
\begin{proof} By Theorem \ref{MonoThm}, we may assume that $(\Tet,O)$ has an obtuse path $A-B-C-D$ and equilibrium on either $\triangle DBC$ or $\triangle ACD$. As the dihedrals on $\overline{BC}$ and $\overline{CD}$ are obtuse, the dihedral on $\overline{BD}$ must be acute. However, by hypothesis, the other two dihedrals at $B$ are obtuse. 

The local geometry at $B$ thus corresponds to a spherical triangle $\triangle \delta\gamma\alpha$ with acute angle at $\delta$, obtuse angles at $\gamma$ and $\alpha$. Then (Lemma \ref{SmallProp}) the edges $\overline{\delta\gamma}$ and $\overline{\delta\alpha}$ are long, and $\overline{\alpha\gamma}$ is short.  Let $E$ be polar to $\bigcirc \delta\alpha$; it lies (Figure \ref{MonoThmFig}) on the great circle polar to $\delta$, which meets $\overline{\delta\gamma}$ at $F$ and $\overline{\delta\alpha}$ at $G$.

But by assumption $(\Tet,O)$ has no equilibrium on $\triangle DBA$, so $\angle\delta\alpha\Omega$ is obtuse; thus $\Omega$ lies on the far side of $\overline{\alpha E}$ and \emph{a fortiori} $\overline{GE}$ from $D$. Thus $\overline{\delta\Omega}$ is long, $\angle DBO$ is obtuse, and $(\Tet,O)$ has no equilibrium on $B$. A similar argument (using the lack of equilibrium on $\triangle ABC$) shows that $(\Tet,O)$ has no equilibrium on $C$.
\end{proof}

Using polar duality, we also get
\begin{tcor}
   If a weighted tetrahedron $(\Tet,O)$ is mono-unstable, it has stable equilibrium on exactly two faces.
\end{tcor}

We can now prove a result analogous to Theorem \ref{MonoThm} for mono-unstable tetrahedra, which does not appear to follow from that result via polar duality. Define an \emph{obtuse cycle} to be a cycle of edges on a tetrahedron $A-B-C-D-A$ such that the face angles $\angle ABC$, $\angle BCD$, and $\angle CDA$ are all obtuse. 

\begin{thm}
    A tetrahedron  $\mathcal{T}$ has an obtuse cycle  if and only if for some $O$ the pair  $(\mathcal{T},O)$ is mono-unstable.
\end{thm}
\begin{proof}
Suppose that $A-B-C-D-A$ is an obtuse cycle, and $P \in \overline{BC}$; then $\angle ABP = \angle ABC$  is obtuse. By the same argument, the angle $\angle DCP$ is obtuse. Moreover, as $\angle ADC$ is obtuse, so is  $\angle ADP$ for $P\in \relint \overline{BC}$ close enough to $C$. $P$ is on the boundary of $\Tet$, but if we let $O$ be an interior point close enough to $P$, angle $\angle ABO$, $\angle DCO$, and $\angle ADO$ will still be obtuse, and $(\Tet,O)$ will have no equilibrium on $B$, $C$, or $D$.

We now show that monoinstability requires the existence of an obtuse cycle. As a triangle has at most one obtuse angle, a tetrahedron has at most four obtuse face angles; and to be mono-unstable it must have one (or more) at each of the three vertices without equilibrium. The obtuse face angles can thus be partitioned among the vertices in only three ways: $\{0,1,1,1\}$, $\{0,1,1,2\}$, and $\{1,1,1,1\}$.

We will represent a vertex with $m$ obtuse face angles and $n$ obtuse dihedrals by the ordered pair $[m,n]$. As every obtuse dihedral has two ends, the sum of $n$  over the vertices is even; and we can only use the pairs $[m,n]$ found in Table 1. The only possibilities for a mono-unstable tetrahedron are: 

\begin{description}
    \item[I] $\{[0,1],[1,1],[1,1],[1,1]\}$;\\
    \item[II] $\{[0,0],[1,1],[1,1],[2,2]\}$;\\
    \item[III] $\{[0,1],[1,1],[1,1],[2,1]\}$;\\
    \item[IV] $\{[0,1],[1,1],[1,1],[2,3]\}$;\\
    \item[V] $\{[1,1],[1,1],[1,1],[1,1]\}$.\\
\end{description}

{\bf I}, which has only three obtuse face angles, is realizable, for instance by a tetrahedron with vertices 
$$\{(-10,0,0),(0,2,0),(0,-2,0),(1,0,1)\}$$ 
(Figure \ref{ObtuseCycle} a).  Let A be the $[0,1]$ vertex, and $\overline{AC}$ its obtuse dihedral, Then $\overline{BD}$ is also an obtuse dihedral, face angles $\angle ABC$, $\angle BCD$, and $\angle CDA$ are obtuse, and $A-B-C-D-A$ is an obtuse cycle.\\

\begin{figure}[h!]
\centering
    \includegraphics[width=0.95\textwidth]{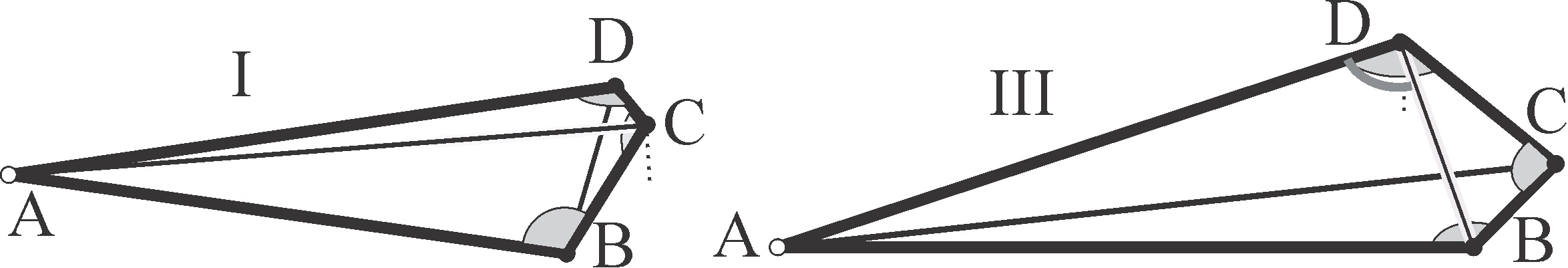}
\caption{Tetrahedra with obtuse cycles}
\label{ObtuseCycle}
\end{figure}

{\bf II} cannot occur.  Let D be the $[2,2]$ vertex, with obtuse angles $\angle ADB$ and $\angle ADC$. Then the dihedrals on $\overline{DC}$ and $\overline{DB}$ are obtuse, B and C are the vertices of type (1,1), and the angles $\angle ABC$ and $\angle ACB$ are both obtuse, which is impossible.\\

{\bf III} can occur. Let A be the $[0,1]$ vertex, D the $[2,1]$ vertex. The tetrahedron has two obtuse dihedrals without a common endpoint. If they were $\overline{AD}$ and $\overline{BC}$, then one of the angles $\angle ADB, \angle ACD$ would be obtuse (without loss of generality $\angle ADB$.) But the angles opposite $\overline{BC}$, that is, $\angle ABD$ and $\angle ACD$, are also obtuse; so $\triangle ABD$ would have two obtuse angles, which is impossible.

However, if (without loss of generality) the obtuse dihedrals are $\overline{AC}$ and $\overline{BD}$, we can construct examples, for instance
$$(A,B,C,D) =((-10,0,0),(2,0,0),(3,2,0),(0,4,1))$$ 
(Figure \ref{ObtuseCycle} b ). Here $A-B-C-D-A$ is the obtuse cycle.

{\bf IV} cannot occur. Let A be the $[2,3]$ vertex; then the dihedrals on $\overline{AB}$, $\overline{AC}$, and $\overline{AD}$ are all obtuse. Without loss of generality let C and D be the vertices of type $[1,1]$; as before, the angles $\angle CDB$ and $\angle BCD$ are both obtuse.\\

{\bf V} would require the tetrahedron to have two disjoint obtuse dihedrals, WLOG  $\overline{AC}$,$\overline{BD}$, opposite the four obtuse angles; but then the skew quadrilateral $\bowtie ABCD $ would have angles summing to more than $2\pi$, which is impossible.    
\end{proof}

\begin{table}\begin{center}
\begin{tabular}{|l|c|c|c|c|c|}
\hline
&\multicolumn{5}{|c|}{Unstable}\\
&\multicolumn{5}{|c|}{equilibria}\\
\hline
&&1&2&3&4 \\
\cline{2-6}
Stable&1&x&$\sqrt{}$&x&x\\
\cline{2-6}
equilibria&2&$\sqrt{}$&$\sqrt{}$&$\sqrt{}$&$\sqrt{}$\\
\cline{2-6}
&3&x&$\sqrt{}$&$\sqrt{}$&$\sqrt{}$\\
\cline{2-6}
&4&x&$\sqrt{}$&$\sqrt{}$&$\sqrt{}$\\
\hline
\end{tabular}\caption{Possible combinations of equilibria}\label{EqTable}
\end{center}\end{table}

\begin{figure}[h!]
\centering
    \includegraphics[width=0.95\textwidth]{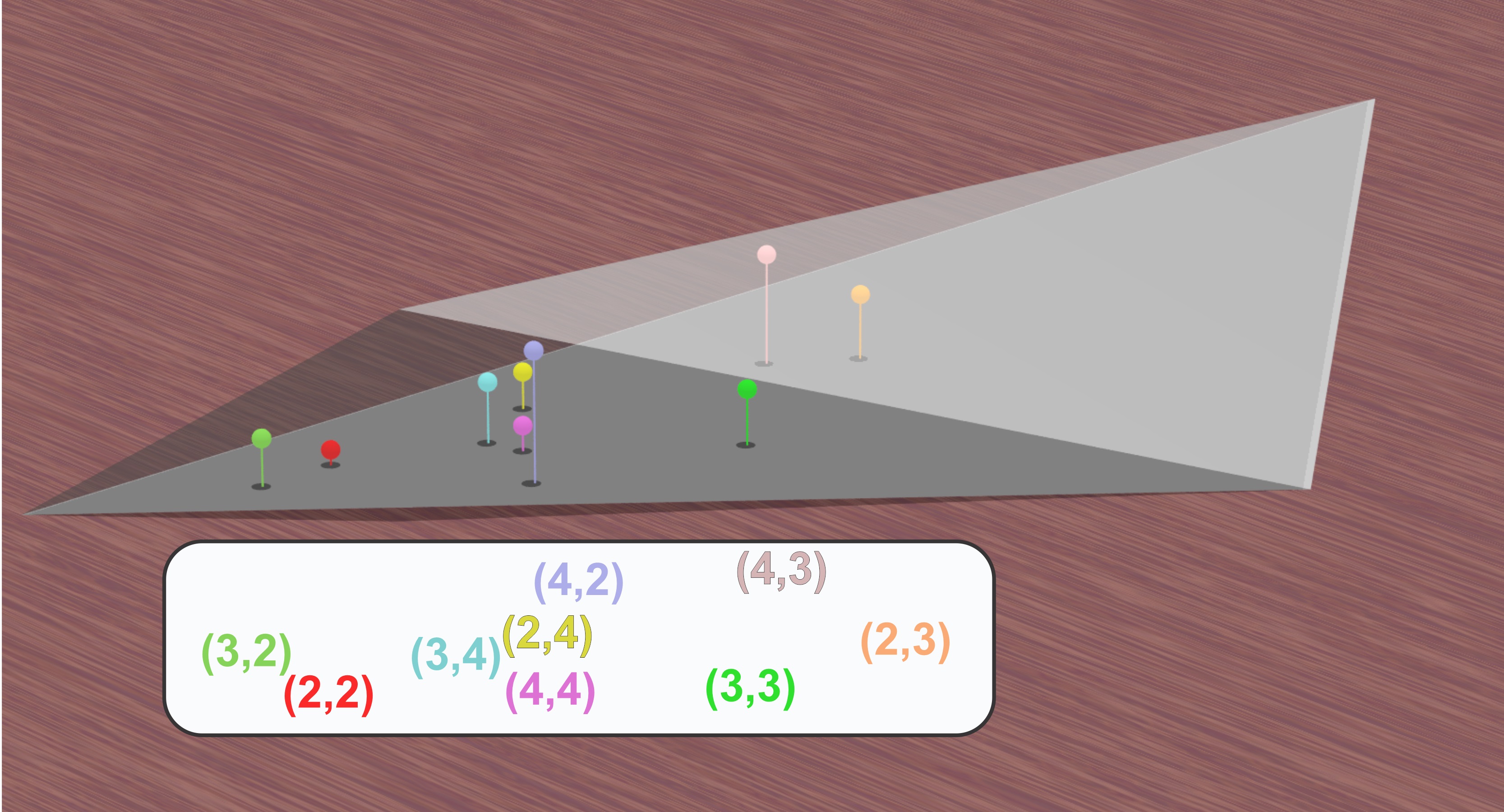}
\caption{A tetrahedron that can have 2--4 stable equilibria and 2--4 unstable equilibria, depending on the choice of center.}
\label{9Centers}
\end{figure}

Examples with every combination of 2--4 stable equilibria and 2--4 unstable equilibria are given in \cite{balancing}. Indeed, there exists a single tetrahedron which exhibits all nine combinations for appropriate choices of centre (Figure \ref{9Centers})\footnote{Vertices are $(0,0,0)$, $(100000,0,0)$, $(50000,41429,0)$, and $(13549,13544,11223)$. Centers are $M_{22}=(15884,5116,835)$, $M_{23}=(46670,11911,3061)$, $M_{24}=(28497,5544,2041)$, $M_{32}=(11400,7243,2597)$, $M_{33}=(33447,17389,3061)$, $M_{34}=(23866,8138,3339)$, $M_{42}=(21845,14097,7142)$, $M_{43}=(42514,9100,6122)$, and $M_{44}=(24407,10239,1391)$.}.

\begin{rem} If a weighted tetrahedron $(\Tet,O)$ is mono-unstable with an equilibrium on vertex $A$, then $A$ has no obtuse face angles. It follows that $(\Tet,O')$ has equilibrium on $A$ for all $O'\in\inter\Tet$; and thus (in contrast to the situation in Theorem 1.3) $\Tet$  cannot be weighted to be mono-unstable on any other vertex. Results like this show that (despite our use of polar duality, and the symmetry of Table \ref{EqTable}) there is no simple duality between stable and unstable equilibria.
\end{rem}

\section{Other polyhedra}
We have seen that no tetrahedron can be mono-monostatic, even when weighted. What about other classes of polyhedron? A vector $(f,e,v)\in\mathbb{N}^3$ is the face vector of some nondegenerate polyhedron if and only if
\begin{itemize} 
    \item $f \geq \frac{v}{2}+2$,
    \item $v \geq \frac{f}{2}+2$,
    \item and $e = f+v-2$. 
\end{itemize}
We'll call such a vector `legal.' We note that equality is obtained in the first expression only when all vertices have degree 3, and the second only when all faces are triangular. 

\begin{thm} Every legal vector except for $(4,6,4)$ is the face vector of a mono-monostatic weighted polyhedron.
\end{thm}
\begin{proof}

Let $\P$ be a weighted polyhedron with at least one nontriangular face. We claim that some vertex $V$ of $\P$ is included in one, two, or three nontriangular faces. For suppose otherwise: intersecting the halfspaces bounded by supporting planes on these faces and containing $P$, we get a convex polyhedron $\Q$ with at least 4 edges at every vertex and at least 4 edges on every face.  Then $e \geq 2v$, $e\geq 2f$, and so for the Euler characteristic we have $\chi(\Q) \leq 0$, an impossibility.

Let $V$ be such a vertex, let $F$ be a nontriangular face including $V$, and let $\delta>0$. We will construct a new polyhedron $\P'$, which shares every vertex of $\P$ except that $V$ is replaced by a new vertex $V'$. Let $G$ be the intersection of the affine hulls of the other nontriangular faces (if any) of $\P$ at $V$: it's an affine subspace of dimension at least 1. Let $H$ be the open halfspace bounded by $\aff F$ that contains $\inter \P$. Then take $V' \in G\cap H\cap B_\delta(V)$: clearly, at least for small $\delta$, $\P'$ has the same number of vertices as $\P$, and one more face.  Moreover, by taking $\delta$ small enough, we can change the orientations of edges and faces by an angle less than any desired $\epsilon >0$. We'll refer to this below as ``face bending.''

Let $(\P,O)$ be a weighted polyhedron; we assume $O$ to be in general position with respect to all edges and face diagonals. For small enough $\delta$, the following are true:
\begin{itemize}
    \item $O\in\inter \P'$
    \item $(\P',O)$ has an equilibrium on a vertex $X$ if and only if $(\P,O)$ has an equilibrium on the corresponding vertex;
    \item $(\P,O)$ has equilibrium on any face other than $F$ if and only if $(\P',O)$ has equilibrium on the corresponding face;
    \item  $(\P,O)$ has equilibrium on $F$ if and only if $(\P',O)$ has equilibrium on $F'$ or $T$ (this requires the foot of the perpendicular from $O$ to $F$ not to lie on the face diagonal that becomes an edge of $F'$);
    \item $(\P',O)$ cannot have equilibrium on both $F'$ and $T$.
\end{itemize}
 
It follows that if there exists a mono-monostatic polyhedron that has  $v$ vertices, and $f$ faces not all triangles, then there exists one with $v$ vertices and $f+1$ faces. As shown in \cite{balancing}, by polar duality there also exists one with $f$ vertices and $v$ faces. 

We conclude the proof using induction. First, we note that there exists a mono-monostatic polyhedron with face vector $(5,8,5)$ (combinatorially equivalent to a square pyramid) - see Figure \ref{MonoMono}.
\begin{figure}[h!]
\centering
    \includegraphics[width=0.95\textwidth]{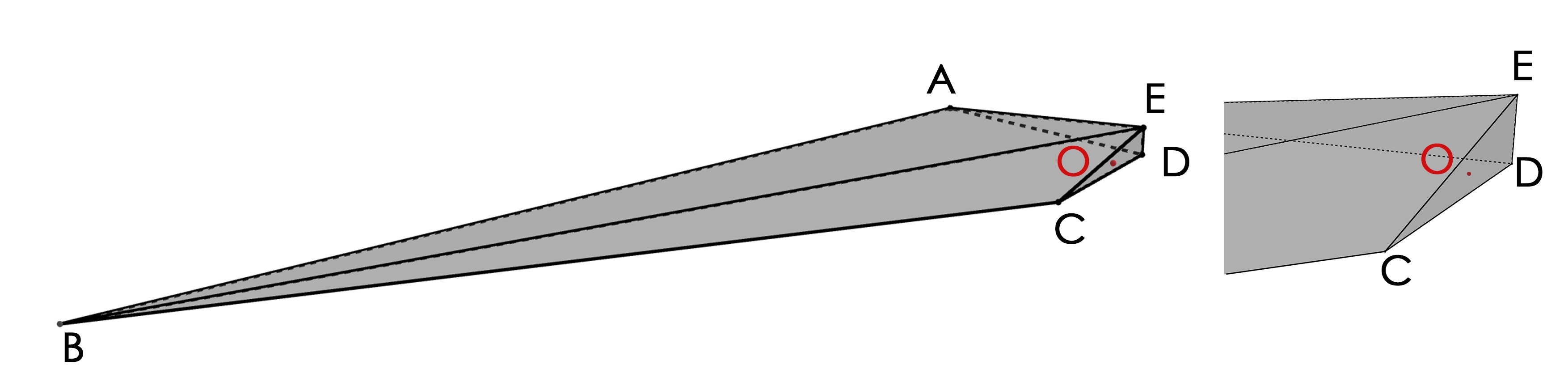}
\caption{A mono-monostatic polyhedron with face vector $(5,8,5)$. Vertices are $(0,0,0)$, $(10000,0,0)$, $(10000,2890,0)$, $(11216,1008,0)$, and $(11216,968,280)$; center of mass is $(10790,643,84)$.}
\label{MonoMono}
\end{figure}

Assume, as an inductive hypothesis, that the claim holds for any legal vector with $v\leq 2n$. Then in particular it holds for $(f,e,v) = (2n+1,3n+2,n+3)$ and $(f,e,v) = (2n+2,3n+3,n+3)$. By polar duality, it also holds for $(f,e,v) = (n+3,3n+2,2n+1)$ and $(f,e,v) = (n+3,3n+3,2n+2)$, both of which minimize $f$ for the given $v$. Face bending then shows that the claim holds for all legal vectors with $v\leq 2n+2$, hence by induction for all legal vectors.
\end{proof}

Except in a few cases, the vector $(f,e,v)$ does not determine the combinatorial class of a polyhedron. We conjecture that, in fact, every combinatorial class of polyhedra, except the tetrahedra, contains elements that admit a mono-monostatic weighting.

\section*{Declaration of competing interests}
The authors have no relevant financial or non-financial interests to disclose.

\end{document}